\documentclass[a4paper,12pt,leqno]{article}
\usepackage{amsmath,amsfonts,amsthm,amssymb,paralist,mathrsfs}
\usepackage[margin=2.5cm]{geometry}

\usepackage[all]{xy}
\numberwithin{equation}{subsection}
\theoremstyle{plain}
\newtheorem{prop}{PROPOSITION}
\newtheorem{corol}{COROLLARY}
\newtheorem{lemma}{LEMMA}
\newtheorem{thm}{THEOREM}
\theoremstyle{definition}
\newtheorem*{defn}{DEFINITION}
\theoremstyle{remark}
\newtheorem*{rem}{REMARK}
\newcommand{\G}{\mathcal{G}}
\newcommand{\F}{\mathcal{F}}
\newcommand{\K}{\mathcal{K}}

\title{K-THEORY OF THE LEAF SPACE OF FOLIATIONS FORMED BY
THE GENERIC K-ORBITS OF SOME INDECOMPOSABLE $MD_5$-GROUPS
}
\author{Le Anh Vu*  and   Duong Quang Hoa**\\
Department of Mathematics and Informatics\\
University of Pedagogy, Ho Chi Minh City, Vietnam\\
\begin{tabular}{rl}
E-mail:&(*) \href{mailto:vula@math.hcmup.edu.vn}{vula@math.hcmup.edu.vn}\\
&(**) \href{mailto:duongquanghoabt@yahoo.com.vn}{duongquanghoabt@yahoo.com.vn}
\end{tabular}}
\date{}
\usepackage{hyperref}
\hypersetup{colorlinks,%
urlcolor=blue,%
linkcolor=black,pdftex}
\DeclareMathOperator{\Ext}{Ext}
\DeclareMathOperator{\Hom}{Hom}
\begin{document}
    \maketitle
    \begin{abstract}
        The paper is a continuation of the authors' work [18]. In [18], we consider foliations formed by the maximal dimensional K-orbits ($MD_5$-foliations) of connected $MD_5$-groups such that their Lie algebras have 4-dimensional commutative derived ideals and give a topological classification of the considered foliations. In this paper, we study K-theory of the leaf space of some of these $MD_5$-foliations and characterize the Connes' C*-algebras of the considered foliations by the method of K-functors.

    \end{abstract}

\footnotetext{

2000AMS Mathematics Subject Classification: Primary 22E45, Secondary
46E25, 20C20.}

    \section*{INTRODUCTION}

    In the decades 1970s $-$ 1980s, works of D.N. Diep [4], J. Rosenberg [10], G. G. Kasparov [7], V. M. Son and H. H. Viet [12],... have seen that K-functors are well adapted to characterize a large class of group C*-algebras. Kirillov's method of orbits allows to find out the class of Lie groups MD, for which the group C*-algebras can be characterized by means of suitable K-functors (see [5]). In terms of D. N. Diep, an MD-group of dimension n (for short, an $MD_n$-group) is an n-dimensional solvable real Lie group whose orbits in the co-adjoint representation (i.e., the K- representation) are the orbits of zero or maximal dimension. The Lie algebra of an $MD_n$-group is called an $MD_n$-algebra (see [5, Section 4.1]).

    In 1982, studying foliated manifolds, A. Connes [3] introduced the notion of C*-algebra associated to a measured foliation. In the case of Reeb foliations (see A. M. Torpe [14]), the method of K-functors has been proved to be very effective in describing the structure of Connes' C*-algebras. For every MD-group G, the family of K-orbits of maximal dimension forms a measured foliation in terms of Connes [3]. This foliation is called MD-foliation associated to G.

    Combining methods of Kirillov (see [8, Section 15]) and Connes (see[3, Section 2, 5]), the first author had studied $MD_4$-foliations associated with all indecomposable connected $MD_4$-groups and characterized Connes' C*-algebras of these foliations in [16]. Recently, Vu and Shum [17] have classified, up to isomorphism, all the $MD_5$-algebras having commutative derived ideals.

    In [18], we have given a topological classification of $MD_5$-foliations associated to the indecomposable connected and simply connected $MD_5$-groups, such that $MD_5$-algebras of them have 4-dimensional commutative derived ideals. There are exactly 3 topological types of the considered $MD_5$-foliations, denoted by $\F_1, \F_2, \F_3$. All  $MD_5$-foliations of type $\F_1$ are the trivial fibrations with connected fibre on 3-dimesional sphere $S^3$, so Connes' C*-algebras of them are isomorphic to the C*-algebra $C(S^3) \otimes \K$ following [3, Section 5], where $\K$ denotes the C*-algebra of compact operators on an (infinite dimensional separable) Hilbert space.

    The purpose of this paper is to study K-theory of the leaf space and to characterize the structure of Connes' C*-algebras $C^*(V,\F)$ of all $MD_5$-foliations $(V,\F)$  of type $\F_2$ by the method of K-functors. Namely, we will express $C^*(V,\F)$  by two repeated extensions of the form

    $$\xymatrix{0 \ar[r] & C_0(X_1)\otimes \K \ar[r] & C^*(V,\F) \ar[r] & B_1 \ar[r] & 0},$$
    $$\xymatrix{0 \ar[r] & C_0(X_2)\otimes \K \ar[r] & B_1 \ar[r] & C_0(Y_2)\otimes \K \ar[r] & 0},$$
then we will compute the invariant system of $C^*(V,\F)$  with
respect to these extensions. If the given C*-algebras are isomorphic
to the reduced crossed products of the form $C_0(V)\rtimes H$, where
H is a Lie group, we can use the Thom-Connes isomorphism to compute
the connecting map $\delta_0, \delta_1$ .

    In another paper, we will study the similar problem for all $MD_5$-foliations of type $\F_3$.

\section{THE $MD_5-$FOLIATIONS OF TYPE $\F_2$}

    Originally, we will recall geometry of K-orbit of $MD_5$-groups which associate with $MD_5$-foliations of type $\F_2$ (see [18]).

    In this section, G will be always an connected and simply connected $MD_5$-group such that its Lie algebras $\G$  is an indecomposable $MD_5$-algebra generated by $\lbrace X_1, X_2, X_3, X_4, X_5 \rbrace$  with $\G^1 : =\left[ {\G,\G} \right] = \mathbb{R}.X_2  \oplus \mathbb{R}.X_3  \oplus \mathbb{R}.X_4  \oplus \mathbb{R}.X_5  \cong \mathbb{R}^4 $, $ad_{X_1} \in End(\G) \equiv Mat_4(\mathbb{R})$. Namely, $\G$  will be one of the following Lie algebras which are studied in [17] and [18].

    1. $\G_{5,4,11 (\lambda_1, \lambda_1 ,\varphi)} $

$$ad_{X_1 }  = \left[ {\begin{array}{*{20}c}
   {\cos \varphi } & { - \sin \varphi } & 0 & 0  \\
   {\sin \varphi } & {\cos \varphi } & 0 & 0  \\
   0 & 0 & {\lambda _1 } & 0  \\
   0 & 0 & 0 & {\lambda _2 }  \\
\end{array}} \right];\lambda _1 ,\lambda _2  \in \mathbb{R} \backslash \left\{ 0 \right\},\lambda _1  \ne \lambda _2 ,\varphi  \in \left( {0,\pi } \right).$$
\newpage
    2. $\G_{5,4,12(\lambda ,\varphi)} $

$$ad_{X_1 }  = \left[ {\begin{array}{*{20}c}
   {\cos \varphi } & { - \sin \varphi } & 0 & 0  \\
   {\sin \varphi } & {\cos \varphi } & 0 & 0  \\
   0 & 0 & \lambda  & 0  \\
   0 & 0 & 0 & \lambda   \\
 \end{array} } \right];\lambda  \in \mathbb{R}\backslash \left\{ 0 \right\},\varphi  \in \left( {0,\pi } \right).$$

    3. $\G_{5,4,13(\lambda ,\varphi)} $

$$ad_{X_1 }  = \left[ {\begin{array}{*{20}c}
   {\cos \varphi } & { - \sin \varphi } & 0 & 0  \\
   {\sin \varphi } & {\cos \varphi } & 0 & 0  \\
   0 & 0 & \lambda  & 1  \\
   0 & 0 & 0 & \lambda   \\
 \end{array} } \right];\lambda  \in \mathbb{R}\backslash \left\{ 0 \right\},\varphi  \in \left( {0,\pi } \right).$$

    The connected and simply connected Lie groups corresponding to these algebras are denoted by $G_{{\text{5,4,11 (}}\lambda _{\text{1}} {\text{,}}\lambda _{\text{1}} {\text{,}}\varphi {\text{) }}} ,{\text{ }}G_{{\text{5,4,12(}}\lambda {\text{,}}\varphi {\text{) }}} ,{\text{ }}G_{{\text{5,4,13(}}\lambda {\text{,}}\varphi {\text{) }}} $. All of these Lie groups are $MD_5$-groups (see [17]) and G is one of  them. We now recall the geometric description of the K-orbits of G in the dual space $\G^*$ of $\G$. Let $\left\{ {X_1^* ,X_2^* ,X_3^* ,X_4^* ,X_5^* } \right\}$   be the basis in $\G^*$ dual to the basis $\left\{ {X_1 ,X_2 ,X_3 ,X_4 ,X_5 } \right\}$ in $\G$. Denote by $\Omega_F$  the K-orbit  of G including $F = \left( {\alpha ,\beta  + i\gamma ,\delta ,\sigma } \right)$  in  $\G^*  \cong \mathbb{R}^5 $.
    \begin{itemize}
        \item If $\beta  + i\gamma  = \delta  = \sigma  = 0$  then $\Omega _F  = \left\{ F \right\}$  (the 0-dimensional orbit).
        \item If  $\left| {\beta  + i\gamma } \right|^2  + \delta ^2  + \sigma ^2  \ne 0$ then $\Omega_F$ is the 2-dimensional orbit as follows

\[
\Omega _F  = \left[ {\begin{array}{*{20}l}
   {\left\{ {\left( {x,\left( {\beta  + i\gamma } \right).e^{\left( {a.e^{ - i\varphi } } \right)} ,\delta .e^{a\lambda _1 } ,\sigma .e^{a\lambda _2 } } \right),{\text{ }}x,a \in \mathbb{R}} \right\}} \hfill  \\
   {{\text{\ \ \ \ \ \ \ \ \ \ \ \ \ \ \ \ \ \ \ \ \ \ \ \ \  when G  =  G}}_{{\text{5,4,11}}\left( {\lambda _{\text{1}} ,\lambda _2 ,\varphi } \right)} ,{\text{ }}\lambda _{\text{1}} ,\lambda _2  \in \mathbb{R}^ *  ,{\text{ }}\varphi  \in \left( {0;\pi } \right).} \hfill  \\
   {\left\{ {\left( {x,\left( {\beta  + i\gamma } \right).e^{\left( {a.e^{ - i\varphi } } \right)} ,\delta .e^{a\lambda } ,\sigma .e^{a\lambda } } \right),{\text{ }}x,a \in \mathbb{R}} \right\}} \hfill  \\
   {{\text{\ \ \ \ \ \ \ \ \ \ \ \ \ \ \ \ \ \ \ \ \ \ \ \ \ when G  =  G}}_{{\text{5,4,12}}\left( {\lambda ,\varphi } \right)} ,{\text{ }}\lambda  \in \mathbb{R}^ *  ,{\text{ }}\varphi  \in \left( {0;\pi } \right).} \hfill  \\
   {\left\{ {\left( {x,\left( {\beta  + i\gamma } \right).e^{\left( {a.e^{ - i\varphi } } \right)} ,\delta .e^{a\lambda } ,\delta .ae^{a\lambda }  + \sigma .e^{a\lambda } } \right),{\text{ }}x,a \in \mathbb{R}} \right\}} \hfill  \\
   {{\text{\ \ \ \ \ \ \ \ \ \ \ \ \ \ \ \ \ \ \ \ \ \ \ \ \  when G  =  G}}_{{\text{5,4,13}}\left( {\lambda ,\varphi } \right)} ,{\text{ }}\lambda  \in \mathbb{R}^ *  ,{\text{ }}\varphi  \in \left( {0;\pi } \right).} \hfill  \\

 \end{array} } \right.
\]

   \end{itemize}

    In [18], we have shown that, the family $\F$ of maximal-dimensional K-orbits of G forms measured foliation in terms of Connes on the open submanifold

    $$V = \left\{ {\left( {x,y,z,t,s} \right) \in G^* :y^2  + z^2  + t^2  + s^2  \ne 0} \right\} \cong \mathbb{R} \times \left( {\mathbb{R}^4 } \right)^* (\subset \G^* \equiv \mathbb{R}^5)$$

    Furthermore, all foliations  $\left( {V,\F_{4,11\left( {\lambda _1 ,\lambda _2 ,\varphi } \right)} } \right)$, $\left( {V,\F_{4,12\left( {\lambda ,\varphi } \right)} } \right)$, $\left( {V,\F_{4,13\left( {\lambda ,\varphi } \right)} } \right)$  are topologically equivalent to each other  $\left( {\lambda _1 ,\lambda _2 ,\lambda  \in \mathbb{R}\backslash \left\{ 0 \right\},\varphi  \in \left( {0;\pi } \right)} \right)$. Thus, we need only choose a envoy among them to describe the structure of the C*-algebra. In this case, we choose the foliation  $\left( {V,\F_{4,12\left( {1,\frac{\pi } {2}} \right)} } \right)$.

    In [18],  we have described the foliation $\left( {V,\F_{4,12\left( {1,\frac{\pi } {2}} \right)} } \right)$  by a suitable action of  $\mathbb{R}^2 $. Namely, we have the following proposition.

    \begin{prop}
        The foliation $\left( {V,\F_{4,12\left( {1,\frac{\pi } {2}} \right)} } \right)$  can be given by an action of the commutative Lie group $\mathbb{R}^2 $  on the manifold $V$.
    \end{prop}

    \begin{proof}

    One needs only to verify that the following action $\lambda $  of $\mathbb{R}^2 $  on V gives the foliation  $\left( {V,\F_{4,12\left( {1,\frac{\pi } {2}} \right)} } \right)$

    \hspace{2.1cm}$\lambda :\mathbb{R}^2  \times V \to V$

    $\left( {\left( {r,a} \right),\left( {x,y + iz,t,s} \right)} \right) \mapsto \left( {x + r,\left( {y + iz} \right).e^{ - ia} ,t.e^a ,s.e^a } \right)$

where  $\left( {r,a} \right) \in \mathbb{R}^2 ,{\text{ }}\left( {x,y + iz,t,s} \right) \in V \cong \mathbb{R} \times \left( {\mathbb{C} \times \mathbb{R}^2 } \right)^ *   \cong \mathbb{R} \times \left( {\mathbb{R}^4 } \right)^ *  $. Hereafter, for simplicity of notation, we write $(V,\F)$ instead of  $\left( {V,\F_{4,12\left( {1,\frac{\pi }{2}} \right)} } \right)$.

    \end{proof}

    It is easy to see that the graph of $\left( {V,\F} \right)$   is indentical with  $V \times \mathbb{R}^2 $, so by [3, Section 5], it follows from Proposition 1 that:

    \begin{corol}[analytical description of $C^*(V,\F)$]
         The Connes’ $C^*$-algebra $ C^*(V,\F)$ can be analytically described the reduced crossed of $C_0 (V)$ by $\mathbb{R}^2 $  as follows

        $$ C^*(V,\F)  \cong C_0 \left( V \right)   \rtimes _\lambda  \mathbb{R}^2.$$ \hfill{$\square$}

    \end{corol}
\hskip 0.5cm
    \section{$C^*(V,\F)$ AS TWO REPEATED EXTENSIONS}
    2.1. Let  $V_1 ,{\text{ }}W_1 ,{\text{ }}V_2 ,{\text{ }}W_2 $ be the following submanifolds of V

    $ V_1  = \left\{ {\left( {x,y,z,t,s} \right) \in V:s \ne 0} \right\} \cong \mathbb{R} \times \mathbb{R}^2  \times \mathbb{R} \times \mathbb{R}^*, $

    $
W_1  = V\backslash V_1  = \left\{ {\left( {x,y,z,t,s} \right) \in V:s = 0} \right\} \cong \mathbb{R} \times \left( {\mathbb{R}^3 } \right)^*  \times \left\{ 0 \right\} \cong \mathbb{R} \times \left( {\mathbb{R}^3 } \right)^*, $

    $
V_2  = \left\{ {\left( {x,y,z,t,0} \right) \in W_1 :t \ne 0} \right\} \cong \mathbb{R} \times \mathbb{R}^2  \times \mathbb{R}^*, $

    $
W_2  = W_1 \backslash V_2  = \left\{ {\left( {x,y,z,t,0} \right) \in W_1 :t = 0} \right\} \cong \mathbb{R} \times \left( {\mathbb{R}^2 } \right)^* .$

    It is easy to see that the action  $\lambda $ in Proposition 1 preserves the subsets  $V_1 ,W_1 , V_2 ,W_2 $. Let $i_1 ,i_2 ,\mu _1 ,\mu _2 $  be the inclusions and the restrictions

    \[
\begin{array}{*{20}c}
   {i_1 :C_0 \left( {V_1 } \right) \to C_0 \left( V \right),} \hfill & {i_2 :C_0 \left( {V_2 } \right) \to C_0 \left( {W_1 } \right),} \hfill  \\
   {\mu _1 :C_0 \left( V \right) \to C_0 \left( {W_1 } \right),} \hfill & {\mu _2 :C_0 \left( {W_1 } \right) \to C_0 \left( {W_2 } \right)} \hfill  \\

 \end{array}
\]
where each function of  $C_0 \left( {V_1 } \right)$  (resp.  $C_0 \left( {V_2 } \right)$) is extented  to the one of  $C_0 \left( V \right)$ (resp.  $C_0 \left( {W_1 } \right)$) by taking the value of zero outside $V_1 $  (resp. $V_2 $).

    It is known a fact that $i_1 ,i_2 ,\mu _1 ,\mu _2 $  are  $\lambda $-equivariant and the following sequences are equivariantly exact:

    \begin{equation}\label{211}
        \xymatrix{0 \ar[r] & C_0(V_1)\ar[r]^{i_1} & C_0(V) \ar[r]^{\mu_1} & C_0(W_1) \ar[r] & 0} \tag{2.1.1}
    \end{equation}
    \begin{equation}
        \xymatrix{0 \ar[r] & C_0(V_2) \ar[r]^{i_2} & C_0(W_1) \ar[r]^{\mu_2} & C_0(W_2) \ar[r] & 0}. \tag{2.1.2}
    \end{equation}

    2.2. Now we denote by  $\left( {V_1,\F_1 } \right),\left( {W_1,\F_1 } \right),\left( {V_2,\F_2 } \right),\left( {W_2,\F_2 } \right)$ the foliations-restrictions of  $(V,\F)$ on $V_1 ,W_1 ,V_2 ,W_2 $  respectively.

    \begin{thm}
    $C^*(V,\F)$ admits the following canonical repeated extensions
        \begin{equation}\label{gamma1}
            \xymatrix{0 \ar[r] & J_1 \ar[r]^{\hspace{-.7cm}\widehat{i_1}} & C^{*}(V,F) \ar[r]^{\hspace{.7cm}\widehat{\mu_1}} & B_1 \ar[r] & 0}, \tag{$\gamma_1$}
        \end{equation}
        \begin{equation}\label{gamma2}
            \xymatrix{0 \ar[r] & J_2 \ar[r]^{\widehat{i_2}} & B_1 \ar[r]^{\widehat{\mu_2}} & B_2\ar[r] & 0}, \tag{$\gamma_2$}
        \end{equation}
where

    $J_1  = C^* \left( {V_1,\F_1 } \right) \cong C_0 \left( {V_1 } \right) \rtimes _\lambda  \mathbb{R}^2  \cong C_0 \left( {\mathbb{R}^3  \cup \mathbb{R}^3 } \right) \otimes K,$

    $J_2  = C^* \left( {V_2,\F_2 } \right) \cong C_0 \left( {V_2 } \right) \rtimes _\lambda  \mathbb{R}^2  \cong C_0 \left( {\mathbb{R}^2  \cup \mathbb{R}^2 } \right) \otimes K,$

    $B_2  = C^* \left( {W_2,\F_2 } \right) \cong C_0 \left( {W_2 } \right) \rtimes _\lambda  \mathbb{R}^2  \cong C_0 \left( {\mathbb{R}_ +  } \right) \otimes K,$

    $B_1  = C^* \left( {W_1,\F_1 } \right) \cong C_0 \left( {W_1 } \right) \rtimes _\lambda  \mathbb{R}^2, $
    and the homomorphismes $\widehat{i_1 },\widehat{i_2 },\widehat{\mu _1 },\widehat{\mu _2 }$  are defined by

        $\left( {\widehat{i_k }f} \right)\left( {r,s} \right) = i_k f\left( {r,s} \right),{\text{      }}k = 1,2$

        $\left( {\widehat{\mu _k }f} \right)\left( {r,s} \right) = \mu _k f\left( {r,s} \right),{\text{      }}k = 1,2$
    \end{thm}

    \begin{proof}
        We note that the graph of $\left( {V_1,\F_1 } \right)$  is indentical with  $V_1  \times \mathbb{R}^2 $, so by [3, section 5],  $J_1  = C^* \left( {V_1,\F_1 } \right) \cong C_0 \left( {V_1 } \right) \rtimes _\lambda  \mathbb{R}^2 $ . Similarly, we have

        $$B_1  \cong C_0 \left( {W_1 } \right) \rtimes _\lambda  \mathbb{R}^2 ,$$

        $$J_2  \cong C_0 \left( {V_2 } \right) \rtimes _\lambda  \mathbb{R}^2 ,$$

        $$B_2  \cong C_0 \left( {W_2 } \right) \rtimes _\lambda  \mathbb{R}^2 ,$$

        From the equivariantly exact sequences in 2.1 and by [2, Lemma 1.1] we obtain the repeated extensions $\left( {\gamma _1 } \right)$  and  $\left( {\gamma _2 } \right)$.

        Furthermore, the foliation $\left( {V_1,\F_1 } \right)$  can be derived from the submersion

        \[
            p_1 :V_1  \approx \mathbb{R} \times \mathbb{R}^2  \times \mathbb{R} \times \mathbb{R}^*  \to \mathbb{R}^3  \cup \mathbb{R}^3 \]
        \[
            p_1 \left( {x,y,z,t,s} \right) = \left( {y,z,t,{\text{sign}}s} \right).
        \]
Hence, by a result of [3, p.562], we get  $J_1  \cong C_0 \left( {\mathbb{R}^3  \cup \mathbb{R}^3 } \right) \otimes K$. The same argument shows that
        \[
J_2  \cong C_0 \left( {\mathbb{R}^2  \cup \mathbb{R}^2 } \right) \otimes K,{\text{          }}B_2  \cong C_0 \left( {\mathbb{R}_ +  } \right) \otimes K.
\]
    \end{proof}

    \section{COMPUTING THE INVARIANT SYSTEM OF $C^*(V,\F)$}

        \begin{defn}
            The set of elements  $\left\{ {\gamma _1 ,\gamma _2 } \right\}$ corresponding to the repeated extensions  $\left( {\gamma _1 } \right)$,  $\left( {\gamma _2 } \right)$ in the Kasparov groups Ext $\left( {B_i ,J_i } \right),{\text{ }}i = 1,2$ is called the system of invariants of $C^* \left( {V,\F} \right)$  and denoted by Index  $C^* \left( {V,\F} \right)$.
        \end{defn}

        \begin{rem}

        Index  $C^* \left( {V,\F} \right)$ determines the so-called “stable type” of  $C^* \left( {V,\F} \right)$ in the set of all repeated extensions
            $$\xymatrix{0\ar[r] & J_1 \ar[r] & E \ar[r] & B_1\ar[r] & 0} ,$$
            $$\xymatrix{0 \ar[r] & J_2 \ar[r] & B_1\ar[r] & B_2 \ar[r] & 0} .$$

        \end{rem}

        The main result of the paper is the following.

        \begin{thm}
            Index  $C^* \left( {V,\F} \right) = \left\{ {\gamma _1 ,\gamma _2 } \right\}$, where

            $\gamma _1  = \left( {\begin{array}{*{20}c}
   0 & 1  \\
   0 & 1  \\
 \end{array} } \right)$ in the group Ext $\left( {B_1 ,J_1 } \right) = Hom\left( {\mathbb{Z}^2 ,\mathbb{Z}^2 } \right)$;

            $\gamma _2  = \left( {1,1} \right)$ in the group Ext $\left( {B_2 ,J_2 } \right) = Hom\left( {\mathbb{Z} ,\mathbb{Z}^2 } \right)$.

        \end{thm}

        To prove this theorem, we need some lemmas as follows.

        \begin{lemma}

            Set $I_2  = C_0 \left( {\mathbb{R}^2  \times \mathbb{R}^* } \right)$ and $A_2  = C_0 \left( {\left( {\mathbb{R}^2 } \right)^* } \right)$

            The following diagram is commutative

            $$\xymatrix{\dots \ar[r] & K_j(I_2)\ar[r] \ar[d]^{\beta_1} & K_j\bigl(C_0(\mathbb{R}^3)^*\bigr) \ar[r] \ar[d]^{\beta_1} & K_j(A_2) \ar[r] \ar[d]^{\beta_1} & K_{j+1}(I_2) \ar[r] \ar[d]^{\beta_1} & \dots\\
            \dots \ar[r] & K_{j+1}\bigl(C_0(V_2)\bigr) \ar[r] & K_{j+1}\bigl(C_0(W_1)\bigr) \ar[r] & K_{j+1}\bigl(C_0(W_2)\bigr)\ar[r] & K_j\bigl(C_0(V_2)\bigr)\ar[r] & \dots}$$
            where  $\beta _1 $ is the isomorphism defined in [13, Theorem 9.7] or in [2, corollary VI.3],  $j \in \mathbb{Z}/2\mathbb{Z}$.
        \end{lemma}
        \begin{proof}
            Let

            $$k_2 :I_2  = C_0 \left( {\mathbb{R}^2  \times \mathbb{R}^* } \right)\xrightarrow{{}}C_0 \left( {\left( {\mathbb{R}^3 } \right)^* } \right)$$
            $$v_2 :C_0 \left( {\left( {\mathbb{R}^3 } \right)^* } \right)\xrightarrow{{}}A_2  = C_0 \left( {\left( {\mathbb{R}^2 } \right)^* } \right)$$
        be the inclusion and restriction defined similarly as in 2.1.

        One gets the exact sequence
            $$\xymatrix{0\ar[r] & I_2 \ar[r]^{\hspace{-.7cm}k_2} & C_0\bigl((\mathbb{R}^3)^*\bigr)\ar[r]^{\hspace{.7cm}v_2} & A_2 \ar[r] & 0}$$

        Note that
        $$C_0 \left( {V_2 } \right) \cong C_0 \left( {\mathbb{R} \times \mathbb{R}^2  \times \mathbb{R}^* } \right) \cong C_0 \left( \mathbb{R} \right) \otimes I_2, $$

        $$C_0 \left( {W_2 } \right) \cong C_0 \left( {\mathbb{R} \times \left( {\mathbb{R}^2 } \right)^* } \right) \cong C_0 \left( \mathbb{R} \right) \otimes A_2, $$

        $$C_0 \left( {W_1} \right) \cong C_0 \left( {\mathbb{R} \times \left( {\mathbb{R}^3 } \right)^* } \right) \cong C_0 \left( \mathbb{R} \right) \otimes C_0 \left( {\mathbb{R}^3 } \right)^* .$$

    The extension (2.1.2) thus can be identified to the following one

            $$\xymatrix{0\ar[r] & C_0(\mathbb{R})\otimes I_2 \ar[r]^{\hspace{-.5cm}id\otimes k_2} & C_0(\mathbb{R})\otimes C_0(\mathbb{R}^3)^* \ar[r]^{\hspace{.5cm}id\otimes v_2} & C_0(\mathbb{R})\otimes A_2\ar[r] & 0}.$$

        Now, using [13, Theorem 9.7; Corollary 9.8] we obtain the assertion of Lemma 1.
        \end{proof}

        \begin{lemma}
        Set $I_1  = C_0 \left( {\mathbb{R}^2  \times \mathbb{R}^* } \right)$ and $A_1  = C\left( {S^2 } \right)$

        The following diagram is commutative

            $$\xymatrix{\dots\ar[r] & K_j(I_1) \ar[r]\ar[d]^{\beta_2} & K_j\bigl(C(S^3)\bigr)\ar[r]\ar[d]^{\beta_2} & K_j(A_1)\ar[r]\ar[d]^{\beta_2} & K_{j+1}(I_1)\ar[r]\ar[d]^{\beta_2} & \dots\\
            \dots \ar[r] & K_j(C_0(V_1))\ar[r] & K_j(C_0(V)) \ar[r] & K_j(C_0(W_1))\ar[r] & K_{j+1}(C_0(V_1))\ar[r] & \dots}$$
        where  $\beta _2 $ is the Bott isomorphism,  $j \in \mathbb{Z}/2\mathbb{Z}$.
        \end{lemma}

        \begin{proof}
            The proof is similar to that of lemma 1, by using the exact sequence (2.1.1) and diffeomorphisms:
            $
V \cong \mathbb{R} \times \left( {\mathbb{R}^4 } \right)^*  \cong \mathbb{R} \times \mathbb{R}_ +   \times S^3 $,
        $
W_1  \cong \mathbb{R} \times \left( {\mathbb{R}^3 } \right)^*  \cong \mathbb{R} \times \mathbb{R}_ +   \times S^2 $.
        \end{proof}

        Before computing the K-groups, we need the following notations. Let $u:\mathbb{R} \to S^1 $ be the map
        $$
            u\left( z \right) = e^{2\pi i\left( {z/\sqrt {1 + z^2 } } \right)} ,{\text{  }}z \in \mathbb{R}
        $$
        Denote by  $u_+  $ (resp.  $u_-$) the restriction of $u$  on  $\mathbb{R}_+$ (resp.  $\mathbb{R}_-$). Note that the class   $\left[ {u_+  } \right]$ (resp.  $\left[ {u_-  } \right]$) is the canonical generator of  $K_1 \left( {C_0 \left( {\mathbb{R}_ +  } \right)} \right) \cong \mathbb{Z}$ (resp.  $K_1 \left( {C_0 \left( {\mathbb{R}_ -  } \right)} \right) \cong \mathbb{Z}$). Let us consider the matrix valued function  $p:\left( {\mathbb{R}^2 } \right)^*  \cong S^1  \times \mathbb{R}_ +   \to M_2 \left( \mathbb{C} \right)$ (resp.  $\overline p :S^2  \cong D/S^1  \to M_2 \left( \mathbb{C} \right)$) defined by:

 $$p\left( {x;y} \right)\left( {resp.{\text{ }}\overline p (x,y)} \right) = \frac{1}
{2}\left( {\begin{array}{*{20}c}
   {1 - \cos \pi \sqrt {x^2  + y^2 } } & {\frac{{x + iy}}
{{\sqrt {x^2  + y^2 } }}\sin \pi \sqrt {x^2  + y^2 } }  \\
   {\frac{{x - iy}}
{{\sqrt {x^2  + y^2 } }}\sin \pi \sqrt {x^2  + y^2 } } & {1 + \cos \pi \sqrt {x^2  + y^2 } }  \\
 \end{array} } \right) .$$

    Then $p$  (resp. $\overline p $) is an idempotent of rank 1 for each  $\left( {x;y} \right) \in {\left( {\mathbb{R}^2 } \right)}^{*}$ (resp.  $\left( {x;y} \right) \in D/S^1 $). Let $\left[ b \right] \in K_0 \left( {C_0 \left( {\mathbb{R}^2 } \right)} \right)$  be the Bott element, [1] be the generator of  $K_0 \left( {C\left( {S^1 } \right)} \right) \cong \mathbb{Z}$.

        \begin{lemma}[See [15, p.234{]}]
            \begin{compactenum}[(i)]
                \item[]
                \item $K_0(B_1)\cong \mathbb{Z}^2,\ K_1(B_1)=0$,
                \item $K_0(J_2)\cong \mathbb{Z}^2$ is generated by $\varphi_0\beta_1\bigl([b]\boxtimes[u_+]\bigr)$ and $\varphi_0\beta_1\bigl([b]\boxtimes[u_-]\bigr);\ K_1(J_2)=0$,
                \item $K_0(B_2)\cong \mathbb{Z}$ is generated by $\varphi_0\beta_1\bigl([1]\boxtimes[u_+]\bigr);\ K_1(B_2)\cong \mathbb{Z}$ is generated by $\varphi_1\beta_1\bigl([p]-[\varepsilon_1]\bigr)$,
            \end{compactenum}
            where $\varphi_j,j\in \mathbb{Z}/2\mathbb{Z}$, is the Thom-Connes isomorphism (see[2]), $\beta_1$ is the isomorphism in Lemma 1, $\varepsilon_1$ is the constant matrix $\biggl(\begin{matrix}
            1 & 0\\
            0 & 0
        \end{matrix}\biggr)$
        and $\boxtimes$ is the external tensor product (see, for example, [2,VI.2]).
        \end{lemma}
        \begin{lemma}
            \begin{compactenum}[(i)]
                \item[]
                \item $K_0\bigl(C^{*}(V,\F)\bigr)\cong\mathbb{Z},\ K_1\bigl(C^*(V,\F)\bigr)\cong\mathbb{Z}$,
                \item $K_0(J_1)=0;\ K_1(J_1)\cong \mathbb{Z}^2$ is generated by $\varphi_1\beta_2\bigl([b]\boxtimes[u_+]\bigr)$ and $\varphi_1\beta_2\bigl([b]\boxtimes[u_-]\bigr)$,
                \item $K_1(B_1)=0;\ K_0(B_1)\cong \mathbb{Z}^2$ is generated by $\varphi_0\beta_2[\bar{1}]$ and $\varphi_0\beta_2\bigl([\bar{p}]-[\varepsilon_1]\bigr)$,
            \end{compactenum}
            where $\bar{1}$ is unit element in $C(S^2)$, $\varphi_0$ is the Thom-Connes isomorphism, $\beta_2$ is the Bott isomorphism.
            \begin{proof}
                \begin{compactenum}[(i)]
                    \item[]
                    \item $K_i\bigl(C^*(V,\F)\bigr)\cong K_i\bigl(C(S^3)\bigr)\cong \mathbb{Z},\ i=0,1$.
                    \item  The proof is similar to (ii) of lemma 3.
                    \item  By [9, p.206], we have
                    $$K_0\bigl(C(S^2)\bigr)=\mathbb{Z}[\bar{1}]+\mathbb{Z}[q],\ \text{where}\ q\in P_2\bigl(C(S^2)\bigr).$$
                    Otherwise, in [9, p.48,53,56]; [13, p.162], one has shown that the map
                    $$dim: K_0\bigl(C(S^2)\bigr)\rightarrow \mathbb{Z}$$
                    is a surjective group homomorphism which satisfied $\dim[\bar{1}]=1,\ \ker(\dim)=\mathbb{Z}$ and non-zero element $q\in P_2\bigl(C(S^2)\bigr)$ in the kernel of the map dim has the form $[q]=[\bar{p}]-[\varepsilon_1]$. Hence, the result is derived straight away because $\beta_2$ and $\varphi_1$ are isomorphisms.
                \end{compactenum}
            \end{proof}
        \end{lemma}

    \textbf{Proof of theorem 2}
    \begin{compactenum}
        \item Computation of $(\gamma_1)$. Recall that the extension $(\gamma_1)$ in theorem 1 gives rise to a six-term exact sequence
        $$\xymatrix{0 = K_0(J_1)\ar[r] & K_0\bigl(C^*(V,F)\bigr)\ar[r] & K_0(B_1)\ar[d]^{\delta_0}\\
        0=K_1(B_1)\ar[u]^{\delta_1} & K_1\bigl(C^*(V,F)\bigr) \ar[l] & K_1(J_1)\ar[l]}$$
        By [11, Theorem 4.14], the isomorphisms
        $$\Ext(B_1,J_1)\cong \Hom\bigl((K_0(B_1),K_1(J_1)\bigr) \cong \Hom(\mathbb{Z}^2,\mathbb{Z}^2)$$
        associates the invariant $\gamma_1\in \Ext(B_1,J_1)$ to the connecting map $\delta_0:K_0(B_1)\rightarrow K_1(J_1)$.

        Since the Thom-Connes isomorphism commutes with $K-$theoretical exact sequence (see[14, Lemma 3.4.3]), we have the following commutative diagram $(j\in \mathbb{Z}/2\mathbb{Z})$:
        $$\xymatrix{\dots \ar[r] & K_j(J_1)\ar[r] & K_j\bigl(C^*(V,F)\bigr)\ar[r] & K_j(B_1)\ar[r]& K_{j+1}(J_1)\ar[r] & \dots\\
        \dots \ar[r] & K_j\bigl(C_0(V_1)\bigr)\ar[r]\ar[u]_{\varphi_j} & K_j\bigl(C_0(V)\bigr)\ar[r]\ar[u]_{\varphi_j} & K_j\bigl(C_0(W_1)\bigr) \ar[r]\ar[u]_{\varphi_j} & K_{j+1}\bigl(C_0(V_1)\bigr)\ar[r]\ar[u]_{\varphi_{j+1}}& \dots}$$
        In view of Lemma 2, the following diagram is commutative
        $$\xymatrix{\dots \ar[r] & K_j\bigl(C_0(V_1)\bigr)\ar[r] & K_1\bigl(C_0(V)\bigr)\ar[r] & K_j\bigl(C_0(W_1)\bigr)\ar[r] & K_{j+1}\bigl(C_1(V_1)\bigr)\ar[r] & \dots\\
        \dots \ar[r] & K_j(I_1)\ar[r]\ar[u]_{\beta_2} & K_j\bigl(C(S^3)\bigr)\ar[r]\ar[u]_{\beta_2} & K_j(A_1)\ar[r]\ar[u]_{\beta_2} & K_{j+1}(I_1)\ar[r]\ar[u]_{\beta_2} & \dots}$$
        Consequently, instead of computing $\delta_0:K_0(B_1)\rightarrow K_1(J_1)$, it is sufficient to compute $\delta_0: K_0(A_1)\rightarrow K_1(I_1)$. Thus, by the proof of Lemma 4, we have to define $\delta_0\bigl([\bar{p}]-[\varepsilon_1]\bigr)=\delta_0\bigl([\bar{p}]\bigr)$ (because $\delta_0\bigl([\varepsilon_1]\bigr)=(0;0)$ and $\delta_0\bigl([\bar{1}]\bigr)=(0;0)$). By the usual definition (see[13, p.170]), for $[\bar{p}]\in K_0(A_1),\ \delta_0\big([\bar{p}]\bigr)=\bigl[e^{2\pi i\tilde{p}}\bigr]\in K_1(I_1)$ where $\tilde{p}$ is a preimage of $\bar{p}$ in (a matrix algebra over) $C(S^3)$, i.e. $v_1\tilde{p}=\bar{p}$.

        We can choose $\tilde{p}(x,y,z)=\dfrac{z}{\sqrt{1+z^2}}\bar{p}(x,y),\ (x,y,z)\in S^3$.

        Let $\tilde{p}_+$ (resp. $\tilde{p}_-$) be the restriction of $\tilde{p}$ on $\mathbb{R}^2\times \mathbb{R}_+$ (resp. $\mathbb{R}^2\times\mathbb{R}_-$). Then we have

        $\delta_0\bigl([\bar{p}]\bigr)=\bigl[e^{2\pi i\tilde{p}}\bigr]=\bigl[e^{2\pi i\tilde{p}_+}\bigr]+\bigl[e^{2\pi i\tilde{p}_-}\bigr]\in K_1\Bigl(C_0\bigl(\mathbb{R}^2\bigr)\otimes C_0\bigl(\mathbb{R}_+\bigr)\Bigr)\oplus K_1\Bigl(C_0\bigl(\mathbb{R}^2\bigr)\otimes C_0\bigl(\mathbb{R}_-\bigr)\Bigr)=K_1(I_1)$

        By [13, Section 4], for each function $f:\mathbb{R}_{\pm}\rightarrow Q_n\widetilde{C_0\bigl(\mathbb{R}^2\bigr)}$ such that $\displaystyle \lim_{x\rightarrow \pm 0}f(t)=\lim_{x\rightarrow \pm\infty}f(t)$, where $Q_n\widetilde{C_0\bigl(\mathbb{R}^2\bigr)}=\Bigl\{a \in M_n\widetilde{C_0\bigl(\mathbb{R}^2\bigr)}, e^{2\pi ia}=Id\Bigr\}$, the class $[f]\in K_1\bigl(C_0(\mathbb{R}^2)\otimes C_0(\mathbb{R}_\pm)\bigr)$ can be determined by $[f]=W_f.[b]\boxtimes[u_\pm]$, where $\displaystyle W_f=\dfrac{1}{2\pi i}\int_{\mathbb{R}_\pm} Tr\bigl(f'(z)f^{-1}(z)\bigr)dz$ is the winding number of $f$.

        By simple computation, we get $\delta_0\bigl([p]\bigr)=[b]\boxtimes[u_+]+[b]\boxtimes [u_-]$. Thus $\gamma_1=\Bigl(\begin{smallmatrix}
        0 & 1\\
        0 & 1
        \end{smallmatrix}\Bigr)\in \Hom_{\mathbb{Z}}(\mathbb{Z}^2,\mathbb{Z}^2)$.
        \item Computation of $(\gamma_2)$. The extension $(\gamma_2)$ gives rise to a six-term exact sequence
        $$\xymatrix{K_0(J_2)\ar[r] & K_0(B_1) \ar[r] & K_0(B_2)\ar[d]^{\delta_0}\\
        K_1(B_2)\ar[u]_{\delta_1} & K_1(B_1)\ar[l] & K_1(J_2)=0\ar[l]}$$
        By [11, Theorem 4.14], $\gamma_2=\delta_1\in \Hom\bigl(K_1(B_2),K_0(J_2)\bigr)=\Hom_{\mathbb{Z}}(\mathbb{Z},\mathbb{Z}^2)$. Similarly to part 1, taking account of Lemma 1 and 3, we have the following commutative diagram $(j\in \mathbb{Z}/2\mathbb{Z})$
        $$\xymatrix{\dots \ar[r] & K_j(J_2)\ar[r] & K_j(B_1)\ar[r] & K_j(B_2)\ar[r] & K_{j+1}(J_2)\ar[r] & \dots\\
        \dots \ar[r] & K_j\bigl(C_0(V_2)\bigr)\ar[r] \ar[u]^{\varphi_j}& K_j\bigl(C_0(W_1)\bigr)\ar[r]\ar[u]^{\varphi_j} & K_j\bigl(C_0(W_2)\bigr)\ar[r]\ar[u]^{\varphi_j} & K_{j+1}\bigl(C_0(V_2)\bigr)\ar[r]\ar[u]^{\varphi_{j+1}} & \dots\\
        \dots \ar[r] & K_{j-1}(I_2)\ar[r]\ar[u]^{\beta_1} & K_{j-1}\bigl(C_0(\mathbb{R}^3)^*\bigr) \ar[r]\ar[u]^{\beta_1} & K_{j-1}(A_2)\ar[r]\ar[u]^{\beta_1} & K_j(I_2)\ar[r]\ar[u]^{\beta_1}& \dots}$$
        Thus we can compute $\delta_0: K_0(A_2)\rightarrow K_1(I_2)$ instead of $\delta_1: K_1(B_2)\rightarrow K_0(J_2)$. By the proof of Lemma 3, we have to define $\delta_0\bigl([p]-[\epsilon_1]\bigr)=\delta_0\bigl([p]\bigr)$ (because $\delta_0\bigl([\epsilon_1]\bigr)=(0,0)$). The same argument as above, we get $\delta_0\bigl([p]\bigr)=[b]\boxtimes[u_+]+[b]\boxtimes[u_-]$. Thus $\gamma_2=(1,1)\in \Hom_{\mathbb{Z}}(\mathbb{Z},\mathbb{Z}^2)\cong \mathbb{Z}^2$. The proof is completed. \hfill {$\square$}

    \end{compactenum}
\hskip 1cm
    \section*{REFERENCES}

\begin{enumerate}
    \item   BROWN, L. G.; DOUGLAS, R. G.; FILLMORE, P. A., Extension of C*-algebra and K-homology, Ann. of Math, 105(1977), 265 - 324.

    \item CONNES, A., An Analogue of the Thom Isomorphism for Crossed Products of a C*-algebra by an Action of $\mathbb{R}$ , Adv. In Math., 39(1981), 31 - 55.

    \item CONNES, A., A Survey of Foliations and Operator Algebras, Proc. Sympos. Pure Mathematics, 38(1982), 521 - 628.

    \item DIEP, D. N., Structure of the group C*-algebra of the group of affine transformations of the line (Russian), Funktsional. Anal. I Prilozhen, 9(1975), 63 - 64.

    \item   DIEP, D. N., Method of Noncommutative Geometry for Group C*-algebras. Reseach Notes in Mathematics Series, Vol.416. Cambridge: Chapman and Hall-CRC Press, 1999.

    \item   KAROUBI, M., K-theory: An introduction, Grund. der Math. Wiss. N0 226, Springer-Verlag, Berlin-Heidelberg-New York, 1978.

    \item KASPAROV, G. G., The operator K-functor and extensions of C*-algebras, Math. USSR Izvestija, 16 (1981), No 3, 513 - 572.

    \item KIRILLOV, A. A., Elements of the Theory of Representations, Springer - Verlag, Berlin - Heidenberg - New York, (1976).

    \item RORDAM, M., LARSEN, F., LAUSTSEN, N., “An Introduction to  $K$-Theory  for C*-Algebras”, Cambridge University Press, United Kingdom, (2000).

    \item ROSENBERG, J., The C*-algebras of some real p-adic solvable groups,  Pacific J. Math, 65 (1976), No 1, 175 - 192.

    \item ROSENBERG, J., Homological invariants of extension of C*-algebras, Proc. Sympos. Pure Math., 38(1982), AMS Providence R.I., 35 - 75.

    \item SON, V. M. ; VIET, H. H., Sur la structure des C*-algebres d’une classe de groupes de Lie, J. Operator Theory, 11 (1984), 77 - 90.

    \item TAYLOR, J. L., "Banach Algebras and Topology", in Algebras in Analysis, pp. 118-186, Academic Press, New York, (1975).

    \item  TORPE, A. M., K-theory for the Leaf Space of Foliations by Reeb Component, J. Func. Anal., 61 (1985), 15-71.

    \item VU, L. A., "On the structure of the  $C^*$-Algebra of the Foliation formed by the  –Orbits of maximal dimendion of the Real Diamond Group", Journal of Operator theory, pp. 227–238 (1990).

    \item VU, L. A., The foliation formed by the K - orbits of Maximal Dimension of the MD4-group, PhD Thesis, Ha Noi (1990) (in Vietnamese).

    \item VU, L. A.; SHUM, K. P., Classification of 5-dimensional MD-algebra having commutative derived ideals, Advances in Algebra and Combinatorics, Singapore: World Scientific, 2008, 353-371.

    \item VU, L. A.; HOA, D. Q., The topology of foliations formed by the generic K-orbits of a subclass of the indecomposable MD5-groups, Science in China, series A: Mathemmatics, Volume 52- Number 2, February 2009, 351-360.
\end{enumerate}
\end{document}